\documentclass[11pt, a4paper]{article}
\usepackage{amsmath,amsfonts,amssymb,bm,hyperref,amsthm}
\numberwithin{equation}{section}
\newtheorem{theorem}{Theorem}[section]
\newtheorem{lemma}[theorem]{Lemma}
\newtheorem{corollary}[theorem]{Corollary}
\theoremstyle{definition}

\newtheorem{remark}[theorem]{Remark}

\begin{document}

\title{Spectral asymmetry of the massless Dirac operator on a 3-torus}
\author{ Robert J.~Downes\thanks{RJD:
Department of Mathematics,
University College London,
Gower Street,
London WC1E~6BT,
UK;
R.Downes@ucl.ac.uk,
\url{http://www.homepages.ucl.ac.uk/\~zcahc37/}
}
\and
Michael Levitin\thanks{ML:
Department of Mathematics,
University of Reading,
Whiteknights,
PO Box 220,
Reading RG6~6AX,
UK;
M.Levitin@reading.ac.uk,
\url{http://www.personal.reading.ac.uk/\~ny901965/}
}
\and
Dmitri Vassiliev\thanks{DV:
Department of Mathematics,
University College London,
Gower Street,
London WC1E~6BT,
UK;
D.Vassiliev@ucl.ac.uk,
\url{http://www.homepages.ucl.ac.uk/\~ucahdva/}
}}

\renewcommand\footnotemark{}

%


\date{1 July 2013}

\maketitle
\begin{abstract}
Consider the massless Dirac operator on a 3-torus equipped with
Euclidean metric and standard spin structure. It is known that
the eigenvalues can be calculated explicitly: the spectrum is
symmetric about zero and zero itself is a double eigenvalue.
The aim of the paper is to develop a perturbation theory for the
eigenvalue with smallest modulus with respect to perturbations
of the metric. Here the application of perturbation techniques
is hindered by the fact that eigenvalues of the massless Dirac
operator have even multiplicity, which is a consequence of this
operator commuting with the antilinear operator of charge
conjugation (a peculiar feature of dimension~3). We derive an
asymptotic formula for the eigenvalue with smallest modulus for
arbitrary perturbations of the metric and present two particular
families of Riemannian metrics for which the eigenvalue with
smallest modulus can be evaluated explicitly. We also establish
a relation between our asymptotic formula and the eta invariant.
\end{abstract}


\section{Introduction}
\label{Introduction}

Let $M$ be a 3-dimensional connected compact oriented manifold
without boundary equipped with a smooth Riemannian metric
$g_{\alpha\beta}$, $\alpha,\beta=1,2,3$ being the tensor
indices. Let $W$ be the corresponding massless Dirac operator,
see Appendix A in \cite{jst_part_b} for definition. There are
two basic examples when the spectrum of $W$ can be calculated
explicitly. The first is the unit torus $\mathbb{T}^3$ equipped
with Euclidean metric. The second is the unit sphere
$\mathbb{S}^3$ equipped with metric induced by the natural
embedding of $\mathbb{S}^3$ in Euclidean space $\mathbb{R}^4$.
In both examples the spectrum turns out to be symmetric about
zero, see Appendix B in \cite{jst_part_b} for details.
Physically, this means that in these two examples there is no
difference between the properties of the particle (massless
neutrino) and antiparticle (massless
antineutrino).

As pointed out in
\cite{atiyah_short_paper,atiyah_part_1,atiyah_part_2,atiyah_part_3},
for a general oriented Riemannian 3-manifold $(M,g)$ there is no
reason for the spectrum of the massless Dirac operator $W$ to be
symmetric. However, producing explicit examples of
spectral asymmetry is a difficult task. To our knowledge, the
only explicit example was constructed in \cite{pfaffle}, with
the example based on the idea of choosing a 3-manifold with flat
metric but highly nontrivial topology. In our paper we take a
different route: we stick with the simplest possible topology
(torus) and create spectral asymmetry
by perturbing the metric.

Further on in this paper we work on the unit torus
$\mathbb{T}^3$ parameterized by cyclic coordinates $x^\alpha$,
$\alpha=1,2,3$, of period $2\pi$.

Suppose first that the metric is Euclidean. Then the massless
Dirac operator corresponding to the standard spin structure
(see formula (A.16) in \cite{jst_part_b}) reads
\begin{equation}
\label{massless Dirac operator Euclidean}
W=
-i
\begin{pmatrix}
\frac\partial{\partial x^3}&\frac\partial{\partial x^1}-i\frac\partial{\partial x^2}
\\
\frac\partial{\partial x^1}+i\frac\partial{\partial x^2}&-\frac\partial{\partial x^3}
\end{pmatrix}.
\end{equation}
The operator (\ref{massless Dirac operator Euclidean}) admits
separation of variables, i.e.~one can seek its eigenfunctions
in the form $v(x)=ue^{im_\alpha x^\alpha}$,  $m\in\mathbb{Z}^3$,
$u\in\mathbb{C}^2$, $u\ne0$,
and calculate the eigenvalues and eigenfunctions explicitly.
The spectrum of the operator
(\ref{massless Dirac operator Euclidean}) is as follows.
\begin{itemize}
\item
Zero is an eigenvalue of multiplicity two.
\item
For each $m\in\mathbb{Z}^3\setminus\{0\}$ we have
the eigenvalue $\|m\|$ and unique (up to rescaling) eigenfunction
of the form $ue^{im_\alpha x^\alpha}$.
\item
For each $m\in\mathbb{Z}^3\setminus\{0\}$ we have
the eigenvalue $\,-\|m\|\,$ and unique (up to rescaling) eigenfunction
of the form $ue^{im_\alpha x^\alpha}$.
\end{itemize}

We now perturb the metric, i.e.~consider a metric
$g_{\alpha\beta}(x;\epsilon)$ the components of which are smooth
functions of coordinates $x^\alpha$,
$\alpha=1,2,3$, and small real parameter $\epsilon$, and which satisfies
\begin{equation}
\label{metric when epsilon is zero}
g_{\alpha\beta}(x;0)=\delta_{\alpha\beta}.
\end{equation}

One way of establishing spectral asymmetry of the perturbed
problem is to compare the asymptotic distribution of large
positive eigenvalues and large negative eigenvalues. As
explained in Section~10 of \cite{jst_part_a}, for a generic
first order differential operator this approach allows one to
establish spectral asymmetry. Unfortunately, the massless Dirac
operator is very special in that the second asymptotic
coefficient of its counting function is zero, see formula (1.23)
in \cite{jst_part_b}, so in the first two approximations in
powers of $\lambda$ its large positive eigenvalues are
distributed the same way as its large negative eigenvalues.
Therefore, in order to demonstrate spectral asymmetry of the
perturbed problem, we will, instead of dealing with large
eigenvalues, deal with small eigenvalues.

\section{Main result}
\label{Main result}

Let $W(\epsilon)$ be the massless Dirac operator corresponding
to the metric $g_{\alpha\beta}(x;\epsilon)$. The difficulty with
applying standard perturbation techniques to the operator
$W(\epsilon)$ is that all its eigenvalues have even multiplicity,
this being a consequence of the fact that
the massless Dirac operator $W(\epsilon)$
commutes with the antilinear operator of charge conjugation
\begin{equation}
\label{antilinear}
v=
\begin{pmatrix}
v_1\\
v_2
\end{pmatrix}
\mapsto
\begin{pmatrix}
-\overline{v_2}\\
\overline{v_1}
\end{pmatrix}
=:\mathrm{C}(v),
\end{equation}
see Property 3 in Appendix A of \cite{jst_part_b}.
In order to overcome this difficulty we develop in Sections
\ref{Perturbation process I: preliminaries}--\ref{Perturbation process III: justification}
a perturbation theory  for the massless Dirac operator which accounts for this
charge conjugation symmetry. We show that perturbation-wise
the double eigenvalues of the massless Dirac operator can be treated as if
they were simple eigenvalues: under perturbation
a double eigenvalue remains\footnote{Here,
of course, it is important that we don't have a magnetic field.
A magnetic field would split up a double eigenvalue, see \cite{solovej}.
The fact that the massless Dirac operator and the charge conjugation operator
do not commute in the presence of a magnetic covector potential
is well known in theoretical physics:
see, for example, formula (2.5) in \cite{hortacsu}.}
a double
eigenvalue and all the usual formulae apply, with only one minor modification.
The minor modification concerns the definition of the
pseudoinverse of the unperturbed operator,
see formulae (\ref{regularized inverse 1})--(\ref{regularized inverse 6}).
Namely, in the definition of the
pseudoinverse
we separate out a two-dimensional eigenspace rather
than a one-dimensional eigenspace.

Given a function $f:\mathbb{T}^3\to\mathbb{C}$, we denote by
\begin{equation}
\label{Fourier transform}
\hat f(m):=
\frac1{(2\pi)^3}
\int_{\mathbb{T}^3}e^{-im_\alpha x^\alpha}f(x)\,dx\,,
\qquad m\in\mathbb{Z}^3,
\end{equation}
its Fourier coefficients. Here  $dx:=dx^1dx^2dx^3$.

Let $\lambda_0(\varepsilon)$ be the eigenvalue of the massless Dirac
operator with smallest modulus and let
\begin{equation}
\label{definition of tensor h}
h_{\alpha\beta}(x):=
\left.
\frac{\partial g_{\alpha\beta}}{\partial\epsilon}
\right|_{\epsilon=0}.
\end{equation}
Further on we raise and lower tensor indices using the Euclidean
metric, which means that raising or lowering a tensor
index doesn't change anything. A repeated tensor index always
indicates summation over the values $1,2,3$.

The following theorem is the main result of our paper.

\begin{theorem}
\label{main theorem}
We have
\begin{equation}
\label{main theorem formula}
\lambda_0(\epsilon)=c\,\epsilon^2+O(\epsilon^3)
\quad\text{as}\quad\epsilon\to0,
\end{equation}
where the constant $c$ is given by the formula
\begin{equation}
\label{formula for c sum}
c=\frac i{16}\,\varepsilon_{\alpha\beta\gamma}
\sum_{m\in\mathbb{Z}^3\setminus\{0\}}
\left(
\delta_{\mu\nu}-\frac{m_\mu m_\nu}{\|m\|^2}
\right)
m_\alpha\,\hat h_{\beta\mu}(m)\,\overline{\hat h_{\gamma\nu}(m)}\,.
\end{equation}
Here $\varepsilon_{\alpha\beta\gamma}$
is the totally antisymmetric quantity, $\varepsilon_{123}:=+1$,
and the overline stands for complex conjugation.
\end{theorem}

Theorem \ref{main theorem} warrants the following remarks.

\begin{itemize}

\item
If the constant $c$ defined by formula (\ref{formula for c sum})
is nonzero, then Theorem \ref{main theorem} tells us that
for sufficiently small nonzero $\epsilon$ the spectrum of our
massless Dirac operator is asymmetric about zero.

\item
Theorem \ref{main theorem} is in agreement with the established
view \cite{hitchin,bar}
that there are no topological obstructions preventing the shift of the zero eigenvalue
of the massless Dirac operator.

\item
Theorem \ref{main theorem} is in agreement with the results of
\cite{ammann2009}. This paper deals with the Dirac operator in
the most general setting. When applied to the case of a compact
oriented Riemannian 3-manifold (not necessarily a 3-torus with
Euclidean metric) with specified spin structure the results of
\cite{ammann2009} tell us that if zero is an eigenvalue of the
Dirac operator, then the metric can be perturbed so that the
zero eigenvalue gets shifted. Furthermore, according to
\cite{ammann2011}, the zero eigenvalue can be shifted by
perturbing the metric on an arbitrarily small open set, which is
also in agreement with our Theorem~\ref{main theorem}.

\item
Put
\[
L_{\gamma\nu\beta\mu}:=\frac{i\varepsilon_{\alpha\beta\gamma}}{(2\pi)^3}
\sum_{m\in\mathbb{Z}^3\setminus\{0\}}
\left(
\delta_{\mu\nu}-\frac{m_\mu m_\nu}{\|m\|^2}
\right)
m_\alpha
\int_{\mathbb{T}^3}
e^{i(x-y)^\alpha m_\alpha}\,(\ \cdot\ )\,dy\,,
\]
\[
P_{\gamma\nu\beta\mu}:=
\frac14
(
L_{\gamma\nu\beta\mu}
+
L_{\nu\gamma\beta\mu}
+
L_{\gamma\nu\mu\beta}
+
L_{\nu\gamma\mu\beta}
).
\]
This gives us a first order pseudodifferential operator $P$
acting in the vector space of rank two symmetric complex-valued
tensor fields, $s_{\beta\mu}\mapsto P_{\gamma\nu\beta\mu}s_{\beta\mu}$.
If we equip this vector space with the natural inner product
$
(r,s):=\int_{\mathbb{T}^3}r_{\alpha\beta}\,\overline{s_{\alpha\beta}}\,dx
$
then it is easy to see that the operator $P$ is formally self-adjoint and formula
(\ref{formula for c sum}) can be rewritten as $c=\frac1{128\pi^3}(Ph,h)$,
where $h$ is defined in accordance with (\ref{definition of tensor h}).
This shows that our coefficient $c$ has a nonlocal (global) nature,
with the source of the non\-locality being the factor
\begin{equation}
\label{source of the nonlocality}
\delta_{\mu\nu}-\frac{m_\mu m_\nu}{\|m\|^2}
\end{equation}
in the symbol of the pseudodifferential operator $P$.
In other words, formula (\ref{formula for c sum})
cannot be rewritten in terms of (linearized) local
differential geometric quantities such as the curvature tensor and
the Cotton tensor.


\item
The rank two tensor (\ref{source of the nonlocality})
can be identified with a linear map in $\mathbb{R}^3$,
$
\ p_\mu\mapsto
\left(
\delta_{\mu\nu}-\frac{m_\mu m_\nu}{\|m\|^2}
\right)
p_\nu
\,$.
This linear map is an orthogonal projection: it projects onto the
plane orthogonal to the covector (momentum) $m$.

\item
Suppose that we are looking at a conformal scaling of
the Euclidean metric,
$g_{\alpha\beta}(x;\epsilon)=e^{2\epsilon\varphi(x)}\delta_{\alpha\beta}$,
where $\varphi:\mathbb{T}^3\to\mathbb{R}$.
Then
$h_{\alpha\beta}(x)=2\varphi(x)\delta_{\alpha\beta}$
and formula (\ref{formula for c sum}) becomes
\begin{equation}
\label{formula for c sum conformal}
c=\frac i{4}\,\varepsilon_{\alpha\beta\gamma}
\sum_{m\in\mathbb{Z}^3\setminus\{0\}}
\left(
\delta_{\mu\nu}-\frac{m_\mu m_\nu}{\|m\|^2}
\right)
m_\alpha\delta_{\beta\mu}\delta_{\gamma\nu}
\left|\hat\varphi(m)\right|^2\,.
\end{equation}
The expression in the RHS of
(\ref{formula for c sum conformal})
is zero because the summand
in $\sum_{m\in\mathbb{Z}^3\setminus\{0\}}$ is odd in $m$.
(Another reason why the expression in the RHS of
(\ref{formula for c sum conformal})
is zero
is that the summand is symmetric in $\beta$,~$\gamma$.)
This agrees with the well-known fact that the zero eigenvalue does not shift
under a conformal scaling of the metric, see Theorem 4.3 in \cite{solovej}.

\item
Suppose that we replace the tensor $h_{\alpha\beta}(x)$
by the tensor $h_{\alpha\beta}(-x)$.
Then $\hat h_{\alpha\beta}(m)$ is replaced by $-\hat h_{\alpha\beta}(-m)$
and, introducing a new summation index $n:=-m$
in formula  (\ref{formula for c sum}),
we see that
the coefficient $c$ changes sign. Physically, this means that
formula  (\ref{formula for c sum}) feels the difference between ``left'' and ``right'',
as one would expect of a formula describing a fermion.

\end{itemize}

The proof of Theorem~ \ref{main theorem}
is given in Section~\ref{Proof of Theorem}.
In Section \ref{Axisymmetric case} we treat the special case
when the metric $g_{\alpha\beta}(x;\epsilon)$ is a function of the coordinate $x^1$ only.
In Sections
\ref{Example of quadratic dependence}
and
\ref{Example of quartic dependence}
we present families of metrics for which
the eigenvalue $\lambda_0(\epsilon)$ can be evaluated explicitly.
Finally, in Section~\ref{The eta invariant} we examine the eta
invariant of our $\epsilon$-dependent massless Dirac operator.

\section{Perturbation process I: preliminaries}
\label{Perturbation process I: preliminaries}

Let $M$ be a 3-dimensional connected compact oriented manifold
without boundary equipped with a smooth Riemannian metric
$g_{\alpha\beta}(x)$, $\alpha,\beta=1,2,3$ being the tensor
indices and $x=(x^1,x^2,x^3)$ being local coordinates. The
perturbation theory developed in this section
and Sections
\ref{Perturbation process II: formal procedure}--\ref{Perturbation process III: justification}
does not assume
that the manifold is necessarily a 3-torus.

We perturb the metric in a smooth manner and denote the
perturbed metric by $g_{\alpha\beta}(x;\epsilon)$, where
$\epsilon$ is a small real parameter. Here we assume that
$g_{\alpha\beta}(x;0)$ is the unperturbed metric described in
the previous paragraph.

By $W_{1/2}(\epsilon)$ we denote the massless Dirac operator on
half-densities corresponding to the metric
$g_{\alpha\beta}(x;\epsilon)$, see Appendix A in
\cite{jst_part_b} for details. We choose to work with the
massless Dirac operator on half-densities $W_{1/2}(\epsilon)$
rather than with the massless Dirac operator $W(\epsilon)$
because we do not want our Hilbert space to depend on
$\epsilon$. The difference between the operators $W(\epsilon)$
and $W_{1/2}(\epsilon)$ is explained in Appendix A of
\cite{jst_part_b}: compare formulae (A.3) and (A.19). The
spectra of the operators $W(\epsilon)$ and $W_{1/2}(\epsilon)$
are the same.

The operator $W_{1/2}(\epsilon)$ is actually not a single operator, but an
equivalence class of operators which differ by the transformation
\begin{equation}
\label{special unitary transformation of Weyl operator}
W_{1/2}(\epsilon)\mapsto R\,W_{1/2}(\epsilon)\, R^*,
\end{equation}
where $R(x;\epsilon)$ is an arbitrary smooth $2\times2$ special
unitary matrix-function.
See Property 4 in Appendix A of \cite{jst_part_b} for a detailed
discussion regarding the transformation
(\ref{special unitary transformation of Weyl operator}),
noting that
the massless Dirac operator on half-densities $W_{1/2}(\epsilon)$
differs from the massless Dirac operator $W(\epsilon)$
only by ``scalar'' factors on the left and on the right --- these
``scalar'' factors commute with matrix-functions $R(x;\epsilon)$ and $R^*(x;\epsilon)$.
Obviously, the transformation
(\ref{special unitary transformation of Weyl operator}) does not
affect the spectrum.
Later on, in Section~\ref{Proof of Theorem},
we will use this gauge degree of freedom
to simplify calculations, see formula (\ref{symmetry of coframe}).

The operator $W_{1/2}(\epsilon)$ acts on 2-columns
$v=\begin{pmatrix}v_1\\ v_2\end{pmatrix}$ of complex-valued
half-densities. Our Hilbert space is $L^2(M;\mathbb{C}^2)$,
which is the vector space of 2-columns of square integrable
half-densities equipped with inner product
\begin{equation}
\label{inner product}
\langle v,w\rangle:=\int_M w^*v\,dx\,.
\end{equation}
The domain of the operator $W_{1/2}(\epsilon)$ is
$H^1(M;\mathbb{C}^2)$, which is the Sobolev space of 2-columns
of half-densities that are square integrable together with
their first partial derivatives. It is known that the operator
$W_{1/2}(\epsilon):H^1(M;\mathbb{C}^2)\to L^2(M;\mathbb{C}^2)$
is self-adjoint and that it has a discrete spectrum, with
eigenvalues accumulating to $+\infty$ and $-\infty$. Note that
here neither the Hilbert space nor the domain depend on
$\epsilon$. It is also known that the eigenfunctions of the operator
$W_{1/2}(\epsilon)$ are infinitely smooth.

The antilinear operator of charge conjugation (\ref{antilinear})
maps any element of $L^2(M;\mathbb{C}^2)$ to an element of
$L^2(M;\mathbb{C}^2)$ and any element of $H^1(M;\mathbb{C}^2)$
to an element of $H^1(M;\mathbb{C}^2)$.
As the massless Dirac operator on half-densities $W_{1/2}(\epsilon)$
differs from the massless Dirac operator $W(\epsilon)$
only by real ``scalar'' factors on the left and on the right,
it also commutes with the operator of charge conjugation:
\begin{equation}
\label{commutes}
\mathrm{C}(W_{1/2}(\epsilon)\,v)=W_{1/2}(\epsilon)\,\mathrm{C}(v)\,,
\end{equation}
$\forall v\in H^1(M;\mathbb{C}^2)$.
Note that the operator of charge conjugation does not itself depend on $\epsilon$.

Observe that formulae (\ref{antilinear}) and (\ref{inner product})
imply the following useful identities:
\begin{equation}
\label{useful identity 1}
\mathrm{C}(\mathrm{C}(v))=-v,
\end{equation}
\begin{equation}
\label{useful identity 2}
\langle v\,,\mathrm{C}(v)\rangle=0,
\end{equation}
\begin{equation}
\label{useful identity 2a}
\langle\mathrm{C}(v)\,,\mathrm{C}(w)\rangle=
\langle w,v\rangle.
\end{equation}

Let
\begin{equation}
\label{asymptotic expansion for massless Dirac on half-densities}
W_{1/2}(\epsilon)=W_{1/2}^{(0)}
+\epsilon W_{1/2}^{(1)}
+\epsilon^2 W_{1/2}^{(2)}
+\ldots
\end{equation}
be the asymptotic expansion of the partial differential operator
$W_{1/2}(\epsilon)$ in powers of the small parameter $\epsilon$.
Obviously, the operators $W_{1/2}^{(k)}$, $k=0,1,2,\ldots$,
are formally self-adjoint first order differential operators
which  commute with the antilinear operator of
charge conjugation (\ref{antilinear}).

Suppose that $\lambda^{(0)}$ is a double eigenvalue of the operator
$W_{1/2}^{(0)}$. As explained in Appendix A of \cite{jst_part_b},
eigenvalues of the massless Dirac operator have even multiplicity,
so a double eigenvalue is the ``simplest'' eigenvalue one can get.

\begin{remark}
The spectrum of the massless Dirac operator on a 3-torus
equipped with Euclidean metric was written down explicitly in
Section~\ref{Introduction}. Examination of the relevant formulae
shows that the only double eigenvalue is the eigenvalue zero as
all others have multiplicity greater than or equal to six.
However,
in this section
and Sections
\ref{Perturbation process II: formal procedure}--\ref{Perturbation process III: justification}
we do not use the fact that
$\lambda^{(0)}=0$.
\end{remark}

Let $v^{(0)}$ be a normalized, $\|v^{(0)}\|=1$, eigenfunction of the
operator $W_{1/2}^{(0)}$ corresponding to the eigenvalue
$\lambda^{(0)}$. Formula (\ref{commutes}) and the fact that
$\lambda^{(0)}$ is real imply that
$\mathrm{C}(v^{(0)})$ is also an eigenfunction
of the
operator $W_{1/2}^{(0)}$ corresponding to the eigenvalue
$\lambda^{(0)}$. Formula (\ref{useful identity 2a}) implies that
$\|\mathrm{C}(v^{(0)})\|=1$,
and, moreover, in view of formula
(\ref{useful identity 2}), the eigenfunctions
$v^{(0)}$
and
$\mathrm{C}(v^{(0)})$
are orthogonal.

The argument presented in the previous paragraph shows that,
when dealing with a double eigenvalue of the massless Dirac operator,
it is sufficient to construct only one eigenfunction: the other one is
obtained by charge conjugation. The argument is valid not only for the
unperturbed operator $W_{1/2}^{(0)}$, but for the perturbed operator
$W_{1/2}(\epsilon)$ as well, provided that $\epsilon$ is small enough
(so that the multiplicity of the eigenvalue does not increase). Hence,
in the perturbation process described in the next section we shall construct
one eigenfunction only.

In the perturbation process that we will describe in the next section
we will make use of the pseudoinverse of the unperturbed operator.
This operator, which we denote by $Q$,
is defined as follows. Consider the problem
\begin{equation}
\label{regularized inverse 1}
(W_{1/2}^{(0)}-\lambda^{(0)})v=f
\end{equation}
where $f\in L^2(M;\mathbb{C}^2)$ is given and
$v\in H^1(M;\mathbb{C}^2)$ is to be found.
Suppose that $f$ satisfies the conditions
\begin{equation}
\label{regularized inverse 2}
\langle f,v^{(0)}\rangle=\langle f,\mathrm{C}(v^{(0)})\rangle=0.
\end{equation}
Then the problem (\ref{regularized inverse 1}) can be resolved
for $v$, however this solution is not unique. We achieve uniqueness
by imposing the conditions
\begin{equation}
\label{regularized inverse 3}
\langle v,v^{(0)}\rangle=\langle v,\mathrm{C}(v^{(0)})\rangle=0
\end{equation}
and define
$Q$ as the linear operator mapping $f$ to $v$,
\begin{equation}
\label{regularized inverse 4}
Q:f\mapsto v.
\end{equation}
Thus,
$Q$ is a bounded linear operator acting in the orthogonal complement of the eigenspace
of the operator
$W_{1/2}^{(0)}$ corresponding to the eigenvalue $\lambda^{(0)}$.
We extend this operator to the whole Hilbert space
$L^2(M;\mathbb{C}^2)$ in accordance with
\begin{equation}
\label{regularized inverse 6}
Qv^{(0)}=Q\mathrm{C}(v^{(0)})=0.
\end{equation}
It is clear from the above definition that the bounded linear operator
$Q$ is self-adjoint and commutes
with the antilinear operator of charge conjugation (\ref{antilinear}).
Note that our definition of the pseudoinverse $Q$ of the unperturbed
operator $W_{1/2}^{(0)}-\lambda^{(0)}$ is in agreement with Rellich's,
see Chapter 2 Section 2 in~\cite{rellich}.

Throughout our perturbation process we will have to deal with various
formally self-adjoint linear operators which commute with
the antilinear operator of charge conjugation (\ref{antilinear}).
Such operators possess a special property which is the subject of
the following lemma.

\begin{lemma}
\label{special property}
Let $L:C^\infty(M;\mathbb{C}^2)\to C^\infty(M;\mathbb{C}^2)$ be a
(possibly unbounded) formally self-adjoint linear operator which commutes with
the antilinear operator of charge conjugation (\ref{antilinear}). Then
for any $v\in C^\infty(M;\mathbb{C}^2)$ we have
\begin{equation}
\label{useful identity 4}
\langle Lv,\mathrm{C}(v)\rangle=0.
\end{equation}
\end{lemma}

\emph{Proof\ }
Take arbitrary $v,w\in C^\infty(M;\mathbb{C}^2)$.
Using formula (\ref{useful identity 2a}) and the fact
that $L$ is formally self-adjoint and commutes with $\mathrm{C}$,
we get
\begin{equation}
\label{useful identity 3}
\langle L\,\mathrm{C}(w)\,,\mathrm{C}(v)\rangle=
\langle\mathrm{C}(L\,w)\,,\mathrm{C}(v)\rangle=
\langle v,Lw\rangle=
\langle Lv,w\rangle.
\end{equation}
For $w=\mathrm{C}(v)$ formula (\ref{useful identity 3}) reads
\begin{equation}
\label{useful identity 3a}
\langle L\,\mathrm{C}(\mathrm{C}(v))\,,\mathrm{C}(v)\rangle=
\langle Lv,\mathrm{C}(v)\rangle.
\end{equation}
But in view of (\ref{useful identity 1}) formula (\ref{useful identity 3a}) can be
rewritten as
\[
-\langle L\,v\,,\mathrm{C}(v)\rangle
=\langle L\,v\,,\mathrm{C}(v)\rangle,
\]
which gives us the required  identity (\ref{useful identity 4}).~$\square$

\section{Perturbation process II: formal procedure}
\label{Perturbation process II: formal procedure}

We now write down the formal perturbation process. A rigorous
justification will be provided in the next section.

Further on in this section as well as in the two following sections
(Sections \ref{Perturbation process III: justification} and \ref{Proof of Theorem})
we write, for the sake of brevity, $A(\epsilon)=W_{1/2}(\epsilon)$
and $A^{(k)}=W_{1/2}^{(k)}$, $k=0,1,2,\ldots$.
In this new notation formula
(\ref{asymptotic expansion for massless Dirac on half-densities}) reads
\begin{equation}
\label{asymptotic expansion for massless Dirac on half-densities brief}
A(\epsilon)=A^{(0)}
+\epsilon A^{(1)}
+\epsilon^2 A^{(2)}
+\ldots.
\end{equation}

We need to solve the eigenvalue problem
\begin{equation}
\label{eigenvalue problem with epsilon}
A(\epsilon)\,v(\epsilon)=\lambda(\epsilon)\,v(\epsilon)\,.
\end{equation}
We seek the eigenvalue and eigenfunction of the perturbed operator
$A(\epsilon)$ in the form of asymptotic expansions
\begin{equation}
\label{asymptotic expansion for eigenvalue}
\lambda(\epsilon)=\lambda^{(0)}
+\epsilon\lambda^{(1)}
+\epsilon^2\lambda^{(2)}
+\ldots,
\end{equation}
\begin{equation}
\label{asymptotic expansion for eigenfunction}
v(\epsilon)=v^{(0)}
+\epsilon v^{(1)}
+\epsilon^2 v^{(2)}
+\ldots.
\end{equation}
Note that we do not aim to preserve the normalization of our eigenfunction
throughout the perturbation process.

Let us forget for a moment that we are dealing with a double eigenvalue
and suppose that our eigenvalue is simple. Then the iterative
procedure for the determination of
$\lambda^{(k)}$ and $v^{(k)}$, $k=1,2,\ldots$, is well known,
see Chapter 2 Section 2 in \cite{rellich}.
At the $k$th step we get the equation
\begin{equation}
\label{step k equation 1}
(A^{(0)}-\lambda^{(0)})v^{(k)}=f^{(k)},
\end{equation}
where
\begin{equation}
\label{step k equation 2}
f^{(k)}:=F^{(k)}v^{(0)},
\end{equation}
and $F^{(k)}$ is some linear operator.
The explicit formula for the operator
$F^{(k)}$ appearing in equations
(\ref{step k equation 1}),
(\ref{step k equation 2})
is written as follows.
Put
\begin{equation}
\label{step k equation 4}
D(\epsilon):=
(B(0)-B(\epsilon))
\Biggl(
I+
\sum_{j=1}^\infty
\bigl[Q\,(B(0)-B(\epsilon))\bigr]^j
\Biggr),
\end{equation}
where $I$ is the identity operator,
$B(\epsilon):=A(\epsilon)-\lambda(\epsilon)I$
and the infinite sum is understood as an asymptotic series.
The operator $D(\epsilon)$ can be expanded in powers of the small parameter~$\epsilon$,
\begin{equation}
\label{step k equation 5}
D(\epsilon)=\sum_{k=1}^\infty\epsilon^kF^{(k)},
\end{equation}
giving us the required $F^{(k)}$. The real number $\lambda^{(k)}$
is determined from the condition
\begin{equation}
\label{step k equation 6}
\langle f^{(k)},v^{(0)}\rangle=0
\end{equation}
after which we resolve (\ref{step k equation 1})
by setting
\begin{equation}
\label{step k equation 7}
v^{(k)}=Qf^{(k)}.
\end{equation}

We claim that the above process carries over to the case of a double
eigenvalue that we are dealing with. Indeed, the difference between
the cases of a simple eigenvalue and a double eigenvalue is that
at the $k$th step of the iterative process in addition to condition
(\ref{step k equation 6}) we need to satisfy the condition
\begin{equation}
\label{step k equation 8}
\langle f^{(k)},\mathrm{C}(v^{(0)})\rangle=0.
\end{equation}
The structure of the operator (\ref{step k equation 4}) is such that
it is formally self-adjoint and commutes with
the antilinear operator of charge conjugation (\ref{antilinear}),
so the operator $F^{(k)}$ defined in accordance with formula
(\ref{step k equation 5})
has the same properties and,
hence, by Lemma~\ref{special property},
condition (\ref{step k equation 8})
is satisfied automatically and the asymptotic process
continues as if the eigenvalue were simple.

We end this section by giving, for future reference, the explicit
formulae for the coefficients
$\lambda^{(1)}$ and $\lambda^{(2)}$
appearing in the asymptotic expansion
(\ref{asymptotic expansion for eigenvalue}):
\begin{equation}
\label{step 1 equation 5}
\lambda^{(1)}=
\langle A^{(1)}v^{(0)},v^{(0)}\rangle,
\end{equation}
\begin{equation}
\label{step 2 equation 7}
\lambda^{(2)}
=
\langle A^{(2)}v^{(0)},v^{(0)}\rangle
-
\langle(A^{(1)}-\lambda^{(1)})\,
Q\,
(A^{(1)}-\lambda^{(1)})\,v^{(0)},v^{(0)}\rangle\,.
\end{equation}

\section{Perturbation process III: justification}
\label{Perturbation process III: justification}

Recall that by $\lambda^{(0)}=\lambda(0)$ we denote a double eigenvalue
of the unperturbed operator $A^{(0)}=A(0)$
(the unperturbed massless Dirac operator on half-densities).
Let us choose a $\delta>0$ such that $\lambda^{(0)}$ is the
only eigenvalue of the operator $A^{(0)}$ on the interval
$[\lambda^{(0)}-\delta,\lambda^{(0)}+\delta]$.

In order to justify our perturbation process we will need the following lemma.

\begin{lemma}
\label{justification}
For sufficiently small $\epsilon$ the interval
\begin{equation}
\label{interval}
(\lambda^{(0)}-\delta,\lambda^{(0)}+\delta)
\end{equation}
contains exactly one double eigenvalue of the operator
$A(\epsilon)$ and no other eigenvalues.
\end{lemma}

\emph{Proof\ }
Denote
$C_\delta:=
\{\mu\in\mathbb{C}\,|\ |\mu-\lambda^{(0)}|=\delta\}$
(circle in the complex plane)
and
$D_\delta:=
\{\mu\in\mathbb{C}\,|\ |\mu-\lambda^{(0)}|<\delta\}$
(open disc in the complex plane).
Put $R_\mu^{(0)}:=(A^{(0)}-\mu I)^{-1}$.
Clearly, for $\mu\in C_\delta$ the operator
$R_\mu^{(0)}$ is well-defined
and, moreover, is a bounded operator acting from
$L^2(M;\mathbb{C}^2)$ to $H^1(M;\mathbb{C}^2)$.
Furthermore, the norm of the operator
$R_\mu^{(0)}:L^2(M;\mathbb{C}^2)\to H^1(M;\mathbb{C}^2)$
is bounded uniformly over $\mu\in C_\delta$.

Let us now define the operator
\begin{equation}
\label{formula for resolvent of the perturbed operator}
R_\mu(\epsilon):=
\left(I+
\sum_{j=1}^\infty
\bigl[-R_\mu^{(0)}(A(\epsilon)-A^{(0)})\bigr]^j
\right)R_\mu^{(0)},
\end{equation}
where $\mu\in C_\delta$.
The operator $A(\epsilon)-A^{(0)}$
is a bounded operator acting from
$H^1(M;\mathbb{C}^2)$ to $L^2(M;\mathbb{C}^2)$
and the norm of the operator
$
A(\epsilon)-A^{(0)}:\linebreak H^1(M;\mathbb{C}^2)\to L^2(M;\mathbb{C}^2)
$
tends to zero as $\epsilon$ tends to zero.
Hence, the series in
(\ref{formula for resolvent of the perturbed operator})
converges for sufficiently small $\epsilon$.
Furthermore, it is easy to see that
\begin{equation}
\label{convergence of resolvents}
R_\mu(\epsilon)\to R_\mu^{(0)}\quad\text{as}\quad\epsilon\to0
\end{equation}
in the sense of the operator norm
$L^2(M;\mathbb{C}^2)\to H^1(M;\mathbb{C}^2)$
and this convergence is uniform over $\mu\in C_\delta$.

Acting onto (\ref{formula for resolvent of the perturbed operator})
with the operator $A(\epsilon)-\mu I$ we see that
$(A(\epsilon)-\mu I)R_\mu(\epsilon)=I$, so
$R_\mu(\epsilon)=(A(\epsilon)-\mu I)^{-1}$.
Put
\begin{equation}
\label{contour integral}
E(\epsilon):=\frac1{2\pi i}\int_{C_\delta}R_\mu(\epsilon)\,d\mu\,.
\end{equation}
The operator $E(\epsilon)$ is the orthogonal projection onto the span
of eigenvectors of the operator $A(\epsilon)$ corresponding to eigenvalues
on the interval (\ref{interval}).
In particular, $E(0)=E^{(0)}$ is the orthogonal projection onto the span
of eigenvectors of the operator $A^{(0)}$ corresponding to the double eigenvalue
$\lambda^{(0)}$.

Formulae (\ref{convergence of resolvents}) and (\ref{contour integral})
imply
\begin{equation}
\label{convergence of spectral projections}
\|E(\epsilon)-E^{(0)}\|_\mathrm{op}\to0\quad\text{as}\quad\epsilon\to0,
\end{equation}
where $\|\,\cdot\,\|_\mathrm{op}$ stands for the operator norm
in the Banach space of bounded linear operators
$L^2(M;\mathbb{C}^2)\to L^2(M;\mathbb{C}^2)$.
Formula (\ref{convergence of spectral projections}) implies that for
sufficiently small $\epsilon$ we have
\begin{equation}
\label{norm is less than one}
\|E(\epsilon)-E^{(0)}\|_\mathrm{op}<1.
\end{equation}
Formula (\ref{norm is less than one})
and the fact that the orthogonal projections
$E(\epsilon)$ and $E^{(0)}$ have finite rank imply that
$\ \operatorname{rank}E(\epsilon)=\operatorname{rank}E^{(0)}=2\,$.
Thus, the operator $A(\epsilon)$ has two eigenvalues, counted with multiplicities,
on the interval (\ref{interval}).
We know, see Property 3 in Appendix A of \cite{jst_part_b},
that the eigenvalues of the operator $A(\epsilon)$ have even multiplicity,
so we are looking at one double eigenvalue on the interval (\ref{interval}).~$\square$

\

Let $\lambda(\epsilon)$ be the unique double eigenvalue
of the operator $A(\epsilon)$ from the interval (\ref{interval}).
Denote by $\sigma(\epsilon)$ the spectrum of the operator $A(\epsilon)$
and, for a given $\mu\in\mathbb{R}$, denote
$\operatorname{dist}(\mu,\sigma(\epsilon))=\min\limits_{\nu\in\sigma(\epsilon)}|\mu-\nu|$.
Obviously, without additional information on $\mu$ and on $\sigma(\epsilon)$
we can only guarantee the inequality
\begin{equation}
\label{spectral inequality}
\operatorname{dist}(\mu,\sigma(\epsilon))\le|\mu-\lambda(\epsilon)|.
\end{equation}

Choose an arbitrary natural $k$ and denote
\begin{equation}
\label{finite asymptotic expansion for eigenvalue}
\tilde\lambda(\epsilon)=\lambda^{(0)}
+\epsilon\lambda^{(1)}
+\epsilon^2\lambda^{(2)}
+\ldots
+\epsilon^k\lambda^{(k)},
\end{equation}
\begin{equation}
\label{finite asymptotic expansion for eigenfunction}
\tilde v(\epsilon)=v^{(0)}
+\epsilon v^{(1)}
+\epsilon^2 v^{(2)}
+\ldots
+\epsilon^k v^{(k)},
\end{equation}
where the $\lambda^{(j)}$ and $v^{(j)}$, $j=0,1,\ldots,k$,
are taken from
(\ref{asymptotic expansion for eigenvalue})
and
(\ref{asymptotic expansion for eigenfunction}).
We have
\begin{equation}
\label{justification 1}
\|\tilde v(\epsilon)\|=1+O(\epsilon),
\end{equation}
\begin{equation}
\label{justification 2}
\|(A(\epsilon)-\tilde\lambda(\epsilon))\tilde v(\epsilon)\|=O(\epsilon^{k+1}),
\end{equation}
where $\|\,\cdot\,\|$ stands for the $L^2(M;\mathbb{C}^2)$ norm
(see (\ref{inner product}) for inner product).
As our operator $A(\epsilon)$ is self-adjoint, formulae
(\ref{justification 1})
and
(\ref{justification 2})
imply
\begin{equation}
\label{justification 3}
\operatorname{dist}(\tilde\lambda(\epsilon),\sigma(\epsilon))
\le
\frac
{\|(A(\epsilon)-\tilde\lambda(\epsilon))\tilde v(\epsilon)\|}
{\|\tilde v(\epsilon)\|}
=O(\epsilon^{k+1}).
\end{equation}

Formulae
(\ref{finite asymptotic expansion for eigenvalue})
and
(\ref{justification 3})
and Lemma~\ref{justification}
imply that for sufficiently small~$\epsilon$,
\begin{equation}
\label{justification 4}
\operatorname{dist}(\tilde\lambda(\epsilon),\sigma(\epsilon))
=
|\tilde\lambda(\epsilon)-\lambda(\epsilon)|,
\end{equation}
compare with (\ref{spectral inequality}).
Combining formulae
(\ref{justification 3})
and
(\ref{justification 4}),
we get
$\lambda(\epsilon)=\tilde\lambda(\epsilon)+O(\epsilon^{k+1})$.
This completes the justification of our perturbation process.

\section{Proof of Theorem \ref{main theorem}}
\label{Proof of Theorem}

The unperturbed massless Dirac operator on half-densities,
which we denote by $A^{(0)}$, is given by the expression in
the RHS of (\ref{massless Dirac operator Euclidean}).
The unperturbed eigenvalue, $\lambda^{(0)}$, is zero
and the corresponding normalized eigenfunction is
\begin{equation}
\label{unperturbed eigenfunction}
v^{(0)}=\frac1{(2\pi)^{3/2}}
\begin{pmatrix}
1\\0
\end{pmatrix}.
\end{equation}
The pseudoinverse $Q$ of the operator $A^{(0)}$ is given by the formula
\begin{multline}
\label{pseudoinverse}
Q
=\frac1{(2\pi)^3}\sum_{m\in\mathbb{Z}^3\setminus\{0\}}
e^{im_\alpha x^\alpha}
\begin{pmatrix}
m_3&m_1-im_2\\
m_1+im_2&-m_3
\end{pmatrix}^{-1}
\int_{\mathbb{T}^3}
e^{-im_\alpha y^\alpha}
\,(\ \cdot\ )\,dy
\\
=\frac1{(2\pi)^3}\sum_{m\in\mathbb{Z}^3\setminus\{0\}}
\frac{e^{im_\alpha x^\alpha}}{\|m\|^2}
\begin{pmatrix}
m_3&m_1-im_2\\
m_1+im_2&-m_3
\end{pmatrix}
\int_{\mathbb{T}^3}
e^{-im_\alpha y^\alpha}
\,(\ \cdot\ )\,dy\,,
\end{multline}
where $dy:=dy^1dy^2dy^3$.
The operator (\ref{pseudoinverse})
is a self-adjoint pseudodifferential operator
of order $-1$.

We have
\begin{equation}
\label{asymptotic expansion for eigenvalue with cubic error}
\lambda(\epsilon)=
\epsilon\lambda^{(1)}
+\epsilon^2\lambda^{(2)}
+O(\epsilon^3),
\end{equation}
where the coefficients $\lambda^{(1)}$ and $\lambda^{(2)}$
are given by formulae
(\ref{step 1 equation 5})
and
(\ref{step 2 equation 7})
respectively. Thus, in order to prove Theorem~ \ref{main theorem}
we need to write down explicitly the differential operators
$A^{(1)}$ and $A^{(2)}$ appearing in the asymptotic expansion
of the perturbed massless Dirac operator on half-densities,
\begin{equation}
\label{asymptotic expansion for operator with cubic error}
A(\epsilon)=
A^{(0)}
+\epsilon A^{(1)}
+\epsilon^2 A^{(2)}
+O(\epsilon^3).
\end{equation}

In what follows we use terminology from microlocal analysis.
In particular, we use the notions of the \emph{principal}
and \emph{subprincipal} symbols of a differential operator,
see subsection 2.1.3 in  \cite{mybook} for details.

Let $L$ be a first order $2\times2$ matrix differential operator.
We denote its principal and subprincipal symbols by
$L_1(x,\xi)$ and $L_\mathrm{sub}(x)$ respectively.
Here $\xi=(\xi_1,\xi_2,\xi_3)$ is the variable dual to the position
variable $x$; in physics literature the $\xi$ would be referred to
as \emph{momentum}.
The subscript in $L_1(x,\xi)$ indicates the degree of homogeneity in $\xi$.

A first order differential operator $L$ is completely determined by its
principal and subprincipal symbols. Indeed, the principal symbol has the form
\begin{equation}
\label{symbols 1}
L_1(x,\xi)=M^{(\alpha)}(x)\,\xi_\alpha\,,
\end{equation}
where $M^{(\alpha)}(x)$ are matrix-functions depending only on the position
variable $x$. It is easy to see that the differential operator $L$ is given by the formula
\begin{equation}
\label{symbols 2}
L=
-\frac i2
M^{(\alpha)}(x)\,\frac\partial{\partial x^\alpha}
-\frac i2
\frac\partial{\partial x^\alpha}\,M^{(\alpha)}(x)
+L_\mathrm{sub}(x)\,.
\end{equation}

Given a first order differential operator $L$, let us consider the expression
$\langle Lv^{(0)},v^{(0)}\rangle$, where $v^{(0)}$ is the constant
column (\ref{unperturbed eigenfunction}) and angular brackets indicate the
inner product (\ref{inner product}). Examination of formula
(\ref{symbols 2}) shows that
\[
\langle Lv^{(0)},v^{(0)}\rangle
=
\langle L_\mathrm{sub}v^{(0)},v^{(0)}\rangle
\]
because the terms coming from the principal symbol integrate to zero.
Consequently
formulae
(\ref{step 1 equation 5})
and
(\ref{step 2 equation 7})
simplify and now read
\begin{equation}
\label{step 1 equation 5 simplified}
\lambda^{(1)}=
\langle A^{(1)}_\mathrm{sub}v^{(0)},v^{(0)}\rangle,
\end{equation}
\begin{equation}
\label{step 2 equation 7 simplified}
\lambda^{(2)}
=
\langle A^{(2)}_\mathrm{sub}v^{(0)},v^{(0)}\rangle
-
\langle(A^{(1)}-\lambda^{(1)})\,
Q\,
(A^{(1)}-\lambda^{(1)})\,v^{(0)},v^{(0)}\rangle\,.
\end{equation}
We see that for the purpose of proving Theorem~ \ref{main theorem}
we do not need to know the full operator
$A^{(2)}$, only its subprincipal symbol $A^{(2)}_\mathrm{sub}$.

In order to write down explicitly the massless Dirac operator on half-densities
$A(\epsilon)$ we need the concepts of \emph{frame} and \emph{coframe}.
The differential geometric definition of coframe was given in Section~3
of \cite{jst_part_b}. However, as in the current paper we are working
in a specified coordinate system, we can adopt a somewhat simpler approach. For the
purposes of the current paper a coframe is a smooth real-valued matrix-function
$e^j{}_\alpha(x;\epsilon)$, $j,\alpha=1,2,3$, satisfying the conditions
\begin{equation}
\label{covariant metric expressed in terms of coframe}
g_{\alpha\beta}(x;\epsilon)
=\delta_{jk}\,e^j{}_\alpha(x;\epsilon)\,e^k{}_\beta(x;\epsilon)\,,
\end{equation}
\begin{equation}
\label{covariant metric expressed in terms of coframe extra}
e^j{}_\alpha(x;0)=\delta^j{}_\alpha\,.
\end{equation}
Here and further on when dealing with matrix-functions we use the convention
that the first index (subscript or superscript) enumerates the rows and
the second index (subscript or superscript) enumerates the columns.
Say, in matrix notation the RHS of
(\ref{covariant metric expressed in terms of coframe})
reads as ``product of coframe transposed and coframe''.

Note that the reason we imposed
condition (\ref{covariant metric expressed in terms of coframe extra})
is so that our unperturbed operator has the form
(\ref{massless Dirac operator Euclidean}).
See also formula
(\ref{special unitary transformation of Weyl operator})
and associated discussion.

For a given metric $g_{\alpha\beta}(x;\epsilon)$ the coframe
$e^j{}_\alpha(x;\epsilon)$ is not defined uniquely. We can
multiply the matrix-function $e^j{}_\alpha(x;\epsilon)$ from the left
by an arbitrary smooth $3\times3$ special orthogonal matrix-function
$O(x;\epsilon)$
satisfying the condition
$O(x;0)=I$,
with $I$ denoting the $3\times3$ identity matrix.
This will
give us a new coframe satisfying the defining relations
(\ref{covariant metric expressed in terms of coframe})
and
(\ref{covariant metric expressed in terms of coframe extra}).
As explained in
Appendix~A of \cite{jst_part_b}, this freedom in the choice of coframe is
a gauge degree of freedom which does not affect the spectrum.
In the current section we specify the gauge by requiring
the matrix-function $e^j{}_\alpha(x;\epsilon)$ to be symmetric,
\begin{equation}
\label{symmetry of coframe}
e^j{}_\alpha(x;\epsilon)=e^\alpha{}_j(x;\epsilon).
\end{equation}
Condition (\ref{symmetry of coframe}) makes sense because
we are working in a specified coordinate system.
Looking ahead, let us point out the main advantage of the symmetric gauge
(\ref{symmetry of coframe}): the asymptotic expansion of the
subprincipal symbol of the massless Dirac operator on half-densities in powers
of $\epsilon$ starts with a quadratic term and, moreover, the coefficient at $\epsilon^2$
has an especially simple structure, see formulae
(\ref{Part 1 of the proof of Theorem eq 1})
and
(\ref{asymptotic expansion for axial torsion with cubic error}).

In matrix notation condition (\ref{covariant metric expressed in terms of coframe}) now reads
``the symmetric positive matrix $g_{\alpha\beta}(x;\epsilon)$ is the square
of the symmetric matrix $e^j{}_\alpha(x;\epsilon)$''.
Conversely, the symmetric matrix $e^j{}_\alpha(x;\epsilon)$
is the square root of the symmetric positive matrix $g_{\alpha\beta}(x;\epsilon)$.
We choose the branch of the square root so that the matrix
$e^j{}_\alpha(x;\epsilon)$ is positive.

According to formulae
(\ref{metric when epsilon is zero})
and
(\ref{definition of tensor h})
we have
\begin{equation}
\label{asymptotic expansion for metric with quadratic error}
g_{\alpha\beta}(x;\epsilon)=
\delta_{\alpha\beta}
+\epsilon h_{\alpha\beta}(x)
+O(\epsilon^2),
\end{equation}
hence, by Taylor's formula for $\sqrt{1+z}\,$,
\begin{equation}
\label{asymptotic expansion for coframe with quadratic error}
e^j{}_\alpha(x;\epsilon)=
\delta^j{}_\alpha
+\frac\epsilon2h^j{}_\alpha(x)
+O(\epsilon^2).
\end{equation}
Here we follow the convention introduced in Section \ref{Main result}:
we raise and lower indices in $h$ using the Euclidean
metric, which means that raising or lowering an
index doesn't change anything. We also
swap, when needed, tensor (Greek) indices
for frame (Latin) indices, which is acceptable because
we are working in a specified coordinate system.

The \emph{frame} is the smooth real-valued matrix-function
$e_j{}^\alpha(x;\epsilon)$, $j,\alpha=1,2,3$, defined by the system
of linear algebraic equations
\begin{equation}
\label{alternative definition of coframe}
e_j{}^\alpha(x;\epsilon)\,e^k{}_\alpha(x;\epsilon)=\delta_j{}^k.
\end{equation}
Note the position of indices in formula (\ref{alternative definition of coframe}).
In matrix notation formula (\ref{alternative definition of coframe}) reads as
``the frame is the transpose of the inverse of the coframe''.
As we chose our coframe to be symmetric, our frame is also symmetric
and is simply the inverse of the coframe.
Formula
(\ref{asymptotic expansion for coframe with quadratic error})
and Taylor's formula for
$(1+z)^{-1}$
imply
\begin{equation}
\label{asymptotic expansion for frame with quadratic error}
e_j{}^\alpha(x;\epsilon)=
\delta_j{}^\alpha
-\frac\epsilon2h_j{}^\alpha(x)
+O(\epsilon^2).
\end{equation}

According to formulae (6.1), (3.5) and (8.1) from \cite{jst_part_b}
the subprincipal symbol of the massless Dirac operator on half-densities is
\begin{equation}
\label{Part 1 of the proof of Theorem eq 1}
A_\mathrm{sub}(x;\epsilon)=
\frac34
\,\bigl({*T}^\mathrm{ax}(x;\epsilon)\bigr)\,I\,,
\end{equation}
where $I$ is the $2\times2$ identity matrix and
${*T}^\mathrm{ax}(x;\epsilon)$ is the scalar function
\begin{multline}
\label{explicit formula for the trace of torsion with a star}
{*T}^\mathrm{ax}
=\frac{\delta_{kl}}3\,\sqrt{\det g^{\alpha\beta}}
\ \bigl[
e^k{}_{1}
\,\partial e^l{}_{3}/\partial x^2
+
e^k{}_{2}
\,\partial e^l{}_{1}/\partial x^3
+
e^k{}_{3}
\,\partial e^l{}_{2}/\partial x^1
\\
-
e^k{}_{1}
\,\partial e^l{}_{2}/\partial x^3
-
e^k{}_{2}
\,\partial e^l{}_{3}/\partial x^1
-
e^k{}_{3}
\,\partial e^l{}_{1}/\partial x^2
\bigr].
\end{multline}
Note that formula
(\ref{asymptotic expansion for metric with quadratic error})
implies
$\,g^{\alpha\beta}(x;\epsilon)=
\delta_{\alpha\beta}
-\epsilon h^{\alpha\beta}(x)
+O(\epsilon^2)\,$,
which, in turn, gives us
\begin{equation}
\label{asymptotic expansion for inverse density with quadratic error}
\sqrt{\det g^{\alpha\beta}(x;\epsilon)}=
1-\frac\epsilon2\operatorname{tr}h(x)
+O(\epsilon^2).
\end{equation}
Substituting formulae
(\ref{asymptotic expansion for inverse density with quadratic error})
and
(\ref{asymptotic expansion for coframe with quadratic error})
into
(\ref{explicit formula for the trace of torsion with a star})
and using the symmetry condition (\ref{symmetry of coframe}),
we get
\begin{multline}
\label{asymptotic expansion for axial torsion with cubic error}
{*T}^\mathrm{ax}(x;\epsilon)=
\\
\frac{\epsilon^2\delta^{kl}}{12}
\left[
h_{k1}
\frac{\partial h_{l3}}{\partial x^2}
+
h_{k2}
\frac{\partial h_{l1}}{\partial x^3}
+
h_{k3}
\frac{\partial h_{l2}}{\partial x^1}
-
h_{k1}
\frac{\partial h_{l2}}{\partial x^3}
-
h_{k2}
\frac{\partial h_{l3}}{\partial x^1}
-
h_{k3}
\frac{\partial h_{l1}}{\partial x^2}
\right]
+O(\epsilon^3)
\\
=
-\frac{\epsilon^2}{12}\varepsilon_{\beta\gamma\delta}
h_{\alpha\beta}\frac{\partial h_{\alpha\gamma}}{\partial x^\delta}
+O(\epsilon^3)\,.
\end{multline}
Formulae
(\ref{Part 1 of the proof of Theorem eq 1})
and
(\ref{asymptotic expansion for axial torsion with cubic error})
imply
\begin{equation}
\label{subprincipal 1}
A^{(1)}_\mathrm{sub}(x)=0,
\end{equation}
\begin{equation}
\label{subprincipal 2}
A^{(2)}_\mathrm{sub}(x)=
-\frac1{16}\,\varepsilon_{\beta\gamma\delta}
\,h_{\alpha\beta}\,\frac{\partial h_{\alpha\gamma}}{\partial x^\delta}\,I\,.
\end{equation}

Substituting
(\ref{subprincipal 1})
into
(\ref{step 1 equation 5 simplified})
we get $\lambda^{(1)}=0$.
Formulae
(\ref{asymptotic expansion for eigenvalue with cubic error})
and
(\ref{step 2 equation 7 simplified})
now simplify and read
\begin{equation}
\label{asymptotic expansion for eigenvalue with cubic error simplified}
\lambda(\epsilon)=
c\,\epsilon^2
+O(\epsilon^3),
\end{equation}
\begin{equation}
\label{step 2 equation 7 simplified further}
c=\lambda^{(2)}
=
\langle A^{(2)}_\mathrm{sub}v^{(0)},v^{(0)}\rangle
-
\langle A^{(1)}\,
Q\,
A^{(1)}\,v^{(0)},v^{(0)}\rangle\,.
\end{equation}
In order to complete our calculation we now need only to write down
the principal symbol $A^{(1)}_1(x,\xi)$ of the differential operator $A^{(1)}$.

According to formulae (A.1)--(A.3) and (A.19) from \cite{jst_part_b}
the principal symbol of the massless Dirac operator on half-densities is
\begin{equation}
\label{principal symbol in terms of frame}
A_1(x,\xi;\epsilon)=
\begin{pmatrix}
e_3{}^\alpha&e_1{}^\alpha-ie_2{}^\alpha
\\
e_1{}^\alpha+ie_2{}^\alpha&-e_3{}^\alpha
\end{pmatrix}
\xi_\alpha\,.
\end{equation}
Formulae
(\ref{principal symbol in terms of frame})
and
(\ref{asymptotic expansion for frame with quadratic error})
imply
\begin{equation}
\label{principal symbol of A1}
A^{(1)}_1(x,\xi)=
-\frac12
\begin{pmatrix}
h_3{}^\alpha&h_1{}^\alpha-ih_2{}^\alpha
\\
h_1{}^\alpha+ih_2{}^\alpha&-h_3{}^\alpha
\end{pmatrix}
\xi_\alpha\,.
\end{equation}
Formulae
(\ref{symbols 1}),
(\ref{symbols 2}),
(\ref{subprincipal 1})
and
(\ref{principal symbol of A1})
allow us to write down the differential operator $A^{(1)}$ explicitly:
\begin{multline}
\label{operator A1}
A^{(1)}=
\frac i4
\begin{pmatrix}
h_3{}^\alpha&h_1{}^\alpha-ih_2{}^\alpha
\\
h_1{}^\alpha+ih_2{}^\alpha&-h_3{}^\alpha
\end{pmatrix}
\frac\partial{\partial x^\alpha}
\\
+
\frac i4
\frac\partial{\partial x^\alpha}
\begin{pmatrix}
h_3{}^\alpha&h_1{}^\alpha-ih_2{}^\alpha
\\
h_1{}^\alpha+ih_2{}^\alpha&-h_3{}^\alpha
\end{pmatrix}
.
\end{multline}

Substituting formulae
(\ref{unperturbed eigenfunction}),
(\ref{pseudoinverse}),
(\ref{subprincipal 2})
and
(\ref{operator A1})
into
(\ref{step 2 equation 7 simplified further})
we arrive at
(\ref{formula for c sum}).
This completes the proof of Theorem~\ref{main theorem}.

\section{Axisymmetric case}
\label{Axisymmetric case}

An important special case is when the metric $g_{\alpha\beta}(x;\epsilon)$
is a function of the coordinate $x^1$ only. In this case one can choose the coframe
and frame so that they depend on the coordinate $x^1$ only and seek eigenfunctions
in the form
\[
v(x)=u(x^1)e^{i(m_2x^2+m_3x^3)},
\qquad m_2,m_3\in\mathbb{Z}.
\]
We get separation of variables, i.e.~the original eigenvalue problem for a partial
differential operator reduces  to
an eigenvalue problem for an ordinary
differential operator depending on the integers $m_2$ and $m_3$ as parameters.
As we know the spectrum of the unperturbed operator
(see Section~\ref{Introduction}) and as the eigenvalues of the original
partial differential operator depend on the small parameter $\epsilon$
continuously, the eigenvalue with smallest modulus will come
from the ordinary differential operator with $m_2=m_3=0$.
We call the case $m_2=m_3=0$ the \emph{axisymmetric case}.

The axisymmetric massless Dirac operator on half-densities reads
\begin{multline}
\label{massless Dirac operator on half-densities axisymmetric}
W_{1/2}(\epsilon)=
\\
-\frac i2
\begin{pmatrix}
e_3{}^1&e_1{}^1-ie_2{}^1
\\
e_1{}^1+ie_2{}^1&-e_3{}^1
\end{pmatrix}
\frac d{d x^1}
-\frac i2
\frac d{d x^1}
\begin{pmatrix}
e_3{}^1&e_1{}^1-ie_2{}^1
\\
e_1{}^1+ie_2{}^1&-e_3{}^1
\end{pmatrix}
\\
+
\frac{\delta_{jk}}{4\sqrt{\det g_{\alpha\beta}}}
\left[
e^j{}_{3}
\left(\frac{d e^k{}_{2}}{d x^1}\right)
-
e^j{}_{2}
\left(\frac{d e^k{}_{3}}{d x^1}\right)
\right]I\,,
\end{multline}
where $I$ is the $2\times2$ identity matrix and
\[
\sqrt{\det g_{\alpha\beta}}
=\frac1{\sqrt{\det g^{\alpha\beta}}}
=\det e^j{}_\alpha
=\frac1{\det e_j{}^\alpha}\,.
\]
Here $e^j{}_\alpha(x^1;\epsilon)$ and $e_j{}^\alpha(x^1;\epsilon)$
are the coframe and frame defined in accordance with formulae
(\ref{covariant metric expressed in terms of coframe}),
(\ref{covariant metric expressed in terms of coframe extra})
and
(\ref{alternative definition of coframe}).

Of course, for a given metric $g_{\alpha\beta}(x^1;\epsilon)$
the coframe $e^j{}_\alpha(x^1;\epsilon)$
and frame $e_j{}^\alpha(x^1;\epsilon)$
are not defined uniquely. We can
multiply the matrix-functions
$e^j{}_\alpha(x^1;\epsilon)$ and $e_j{}^\alpha(x^1;\epsilon)$
from the left
by an arbitrary smooth $3\times3$ special orthogonal matrix-function
$O(x^1;\epsilon)$
satisfying the condition
$O(x^1;0)=I$,
with $I$ denoting the $3\times3$ identity matrix.
This will
give us a new coframe and a new frame satisfying the defining relations
(\ref{covariant metric expressed in terms of coframe}),
(\ref{covariant metric expressed in terms of coframe extra})
and
(\ref{alternative definition of coframe}).
Note that in writing down formula
(\ref{massless Dirac operator on half-densities axisymmetric})
we did not assume a particular choice of gauge,
compare with (\ref{symmetry of coframe}).


In the axisymmetric case formula
(\ref{formula for c sum})
also simplifies and reads now
\begin{equation}
\label{formula for c sum axisymmetric}
c=-\frac18
\sum_{m_1\in\mathbb{N}}
m_1
\operatorname{tr}
\left[
\begin{pmatrix}
\hat h_{22}&\hat h_{23}\\
\hat h_{32}&\hat h_{33}
\end{pmatrix}
\begin{pmatrix}
0&-i\\
i&0
\end{pmatrix}
\begin{pmatrix}
\hat h_{22}&\hat h_{23}\\
\hat h_{32}&\hat h_{33}
\end{pmatrix}^*
\,\right],
\end{equation}
where $\hat h_{\alpha\beta}=\hat h_{\alpha\beta}(m_1)$ and the star
stands for Hermitian conjugation.


\section{Example of quadratic dependence on~$\epsilon$}
\label{Example of quadratic dependence}

Consider the metric
\begin{multline}
\label{metric quadratic}
g_{\alpha\beta}(x^1;\epsilon)\,dx^\alpha dx^\beta
=\bigl[dx^1\bigr]^2
+\bigl[
\bigl(1+\epsilon\bigl(\cos x^1\bigr)\bigr)dx^2
+\epsilon\bigl(\sin x^1\bigr)dx^3
\bigr]^2
\\
+\bigl[
\epsilon\bigl(\sin x^1\bigr)dx^2
+\bigl(1-\epsilon\bigl(\cos x^1\bigr)\bigr)dx^3
\bigr]^2.
\end{multline}
Then
\begin{equation}
\label{Example of quadratic dependence coframe}
e^j{}_\alpha(x^1;\epsilon)
=
\delta^j{}_\alpha
+\epsilon
\begin{pmatrix}
0&0&0\\
0&\cos x^1&\sin x^1\\
0&\sin x^1&-\cos x^1
\end{pmatrix}
\end{equation}
is a coframe associated with the metric
(\ref{metric quadratic}), see formulae
(\ref{covariant metric expressed in terms of coframe}),
(\ref{covariant metric expressed in terms of coframe extra}),
and
\begin{equation}
\label{Example of quadratic dependence frame}
e_j{}^\alpha(x^1;\epsilon)
=
\begin{pmatrix}
1&0&0\\
0&\frac{1-\epsilon\cos x^1}{1-\epsilon^2}&-\frac{\epsilon\sin x^1}{1-\epsilon^2}\\
0&-\frac{\epsilon\sin x^1}{1-\epsilon^2}&\frac{1+\epsilon\cos x^1}{1-\epsilon^2}
\end{pmatrix}
\end{equation}
is the corresponding frame, see formula
(\ref{alternative definition of coframe}).
Note that in writing formula
(\ref{Example of quadratic dependence frame})
we used the fact that
\begin{equation}
\label{Example of quadratic dependence density}
\det e^j{}_\alpha(x;\epsilon)=1-\epsilon^2=\sqrt{\det g_{\alpha\beta}}\,.
\end{equation}
Substituting formulae
(\ref{Example of quadratic dependence coframe})--(\ref{Example of quadratic dependence density})
into (\ref{massless Dirac operator on half-densities axisymmetric})
we get
\begin{equation}
\label{Example of quadratic dependence operator}
W(\epsilon)=
-i
\begin{pmatrix}
0&1\\
1&0
\end{pmatrix}
\frac d{dx^1}
-\frac{\epsilon^2}{2(1-\epsilon^2)}I.
\end{equation}
Note that in the LHS we dropped the subscript $1/2\,$: as the Riemannian density
is constant, see
(\ref{Example of quadratic dependence density}),
there is no need to distinguish the massless
Dirac operator $W(\epsilon)$ and the massless Dirac operator on
half-densities $W_{1/2}(\epsilon)$.

It is easy to see that the eigenvalues of the ordinary differential operator
(\ref{Example of quadratic dependence operator})
subject to the (boundary) condition of $2\pi$-periodicity
are
\[
\lambda_n(\epsilon)=n-\frac{\epsilon^2}{2(1-\epsilon^2)},
\qquad n\in\mathbb{Z},
\]
and that all eigenvalues have multiplicity two. In particular,
the eigenvalue with smallest modulus is
\begin{equation}
\label{eigenvalue quadratic asymptotic formula}
\lambda_0(\epsilon)
=-\frac{\epsilon^2}{2(1-\epsilon^2)}
=-\frac{\epsilon^2}{2}+O(\varepsilon^4)
\quad\text{as}\quad\varepsilon\to0.
\end{equation}

Let us now test Theorem~\ref{main theorem}
by comparing the asymptotic formula from this theorem
with formula (\ref{eigenvalue quadratic asymptotic formula}).
Substituting (\ref{metric quadratic}) into (\ref{definition of tensor h}) we get
\[
h_{\alpha\beta}(x^1)=2
\begin{pmatrix}
0&0&0\\
0&\cos x^1&\sin x^1\\
0&\sin x^1&-\cos x^1
\end{pmatrix}.
\]
Application of the Fourier transform (\ref{Fourier transform}) gives us
\begin{equation}
\label{Example of quadratic dependence h hat}
\hat h_{\alpha\beta}(m_1)=
\begin{cases}
\begin{pmatrix}
0&0&0\\
0&1&-i\\
0&-i&-1
\end{pmatrix}
\quad&\text{for}\quad m_1=1,
\\
{}&{}
\\
\qquad\quad
0\quad&\text{for}\quad m_1=2,3,\ldots.
\end{cases}
\end{equation}
Substituting
(\ref{Example of quadratic dependence h hat})
into
(\ref{formula for c sum axisymmetric})
we get
$c=-\frac12$, in agreement with  (\ref{eigenvalue quadratic asymptotic formula}).

\section{Example of quartic dependence on~$\epsilon$}
\label{Example of quartic dependence}

Consider the metric
\begin{multline}
\label{metric quartic}
g_{\alpha\beta}(x^1;\epsilon)\,dx^\alpha dx^\beta
=\bigl[dx^1+\epsilon\bigl(\cos x^1\bigr)dx^2+\epsilon\bigl(\sin x^1\bigr)dx^3\bigr]^2
\\
+\bigl[dx^2\bigr]^2+\bigl[dx^3\bigr]^2.
\end{multline}
Then
\begin{equation}
\label{Example of quartic dependence coframe}
e^j{}_\alpha(x^1;\epsilon)
=
\delta^j{}_\alpha
+\epsilon
\begin{pmatrix}
0&\cos x^1&\sin x^1\\
0&0&0\\
0&0&0
\end{pmatrix}
\end{equation}
is a coframe associated with the metric
(\ref{metric quartic}), see formulae
(\ref{covariant metric expressed in terms of coframe}),
(\ref{covariant metric expressed in terms of coframe extra}),
and
\begin{equation}
\label{Example of quartic dependence frame}
e_j{}^\alpha(x^1;\epsilon)
=
\delta_j{}^\alpha
-\epsilon
\begin{pmatrix}
0&0&0\\
\cos x^1&0&0\\
\sin x^1&0&0\\
\end{pmatrix}
\end{equation}
is the corresponding frame, see formula
(\ref{alternative definition of coframe}).
Note that in writing formula
(\ref{Example of quartic dependence frame})
we used the fact that
\begin{equation}
\label{Example of quartic dependence density}
\det e^j{}_\alpha(x;\epsilon)=1=\sqrt{\det g_{\alpha\beta}}\,.
\end{equation}
Substituting formulae
(\ref{Example of quartic dependence coframe})--(\ref{Example of quartic dependence density})
into (\ref{massless Dirac operator on half-densities axisymmetric})
we get
\begin{multline}
\label{Example of quartic dependence operator}
W(\epsilon)=
-i
\begin{pmatrix}
0&1\\
1&0
\end{pmatrix}
\frac d{dx^1}
-\frac{\epsilon^2}4I
\\
+\frac{i\epsilon}2
\begin{pmatrix}
\sin x^1&-i\cos x^1
\\
i\cos x^1&-\sin x^1
\end{pmatrix}
\frac d{d x^1}
+\frac{i\epsilon}2
\frac d{d x^1}
\begin{pmatrix}
\sin x^1&-i\cos x^1
\\
i\cos x^1&-\sin x^1
\end{pmatrix}.
\end{multline}
Note that in the LHS we dropped the subscript $1/2\,$: as the Riemannian density
is constant, see
(\ref{Example of quartic dependence density}),
there is no need to distinguish the massless
Dirac operator $W(\epsilon)$ and the massless Dirac operator on
half-densities $W_{1/2}(\epsilon)$.

We shall now rewrite the ordinary differential operator
(\ref{Example of quartic dependence operator})
in a somewhat more convenient form. To this end, let us
introduce the special unitary matrix
\begin{equation}
\label{matrix R}
R:=
\frac1{\sqrt2}
\begin{pmatrix}
1&1
\\
-1&1
\end{pmatrix}
\end{equation}
and put
\begin{equation}
\label{definition of operator W with a tilde}
\tilde W(\epsilon):=R\ W(\epsilon)\,R^*\,,
\end{equation}
compare with formula
(\ref{special unitary transformation of Weyl operator}).
Clearly, the operator $\tilde W(\epsilon)$ has the same
spectrum as the operator $W(\epsilon)$.
Substituting
(\ref{Example of quartic dependence operator})
and
(\ref{matrix R})
into
(\ref{definition of operator W with a tilde})
we arrive at the following explicit formula for
the ordinary differential operator $\tilde W(\epsilon)$:
\begin{multline}
\label{Example of quartic dependence operator with tilde}
\tilde W(\epsilon)=
-i
\begin{pmatrix}
1&0\\
0&-1
\end{pmatrix}
\frac d{dx^1}
-\frac{\epsilon^2}4I
\\
+\frac{i\epsilon}2
\begin{pmatrix}
0&-ie^{-ix^1}
\\
ie^{ix^1}&0
\end{pmatrix}
\frac d{d x^1}
+\frac{i\epsilon}2
\frac d{d x^1}
\begin{pmatrix}
0&-ie^{-ix^1}
\\
ie^{ix^1}&0
\end{pmatrix}.
\end{multline}

The coefficients of the ordinary differential operator
(\ref{Example of quartic dependence operator with tilde})
are trigonometric polynomials and one would not normally expect
the eigenfunctions to be trigonometric polynomials. However,
the operator
(\ref{Example of quartic dependence operator with tilde})
has a special structure which ensures that the eigenfunctions
are trigonometric polynomials. Namely, put
\begin{equation}
\label{eigenvalues of the perturbed operator with tilde}
\lambda_n(\epsilon)=
-\frac12-\frac{\epsilon^2}4
+\sqrt{1+\epsilon^2}\left(n+\frac12\right),
\qquad n\in\mathbb{Z},
\end{equation}
\begin{equation}
\label{eigenfunctions of the perturbed operator with tilde}
v^{(n)}(x^1;\epsilon)
=
\begin{pmatrix}
\left(1+\sqrt{1+\epsilon^2}\,\right)e^{inx^1}
\\
-i\,\epsilon\,e^{i(n+1)x^1}
\end{pmatrix},
\qquad n\in\mathbb{Z}.
\end{equation}
It is easy to see that the column-functions
(\ref{eigenfunctions of the perturbed operator with tilde})
are eigenfunctions of the operator
(\ref{Example of quartic dependence operator with tilde})
corresponding to eigenvalues
(\ref{eigenvalues of the perturbed operator with tilde}).
Moreover, it is easy to see that the charge conjugates,
$\mathrm{C}(v^{(n)}(x^1;\epsilon))$, of the
column-functions
(\ref{eigenfunctions of the perturbed operator with tilde})
are eigenfunctions of the operator
(\ref{Example of quartic dependence operator with tilde})
corresponding to the same eigenvalues
(\ref{eigenvalues of the perturbed operator with tilde}).
This means that the numbers
(\ref{eigenvalues of the perturbed operator with tilde})
are eigenvalues of the operator
(\ref{Example of quartic dependence operator with tilde})
of multiplicity at least two.
Finally, it is easy to see that
\begin{multline}
\label{spans are equal}
\operatorname{span}
\bigl\{
v^{(n)}(x^1;\epsilon),
\ \mathrm{C}(v^{(n)}(x^1;\epsilon))
\ \bigr|
\ n\in\mathbb{Z}
\bigr\}
\\
=
\operatorname{span}
\bigl\{
v^{(n)}(x^1;0),
\ \mathrm{C}(v^{(n)}(x^1;0))
\ \bigr|
\ n\in\mathbb{Z}
\bigr\},
\end{multline}
where $\,\operatorname{span}S\,$ denotes the linear span,
i.e.~set of all finite linear combinations of elements
of a given set $S$. Formula (\ref{spans are equal})
implies that we haven't missed any eigenvalues,
that is, that the list
(\ref{eigenvalues of the perturbed operator with tilde})
contains \emph{all} the eigenvalues of the operator
(\ref{Example of quartic dependence operator with tilde})
and that each of these eigenvalues
has multiplicity two.

\begin{remark}
We do not fully understand the underlying reasons why the
axi\-symmetric massless
Dirac operator corresponding to the metric
(\ref{metric quartic})
admits an explicit evaluation of the eigenvalues and eigenfunctions.
Somehow, this particular Dirac operator has properties similar to those
of an integrable system.
\end{remark}

The eigenvalue
(\ref{eigenvalues of the perturbed operator with tilde})
with smallest modulus is
\begin{equation}
\label{eigenvalue quartic asymptotic formula}
\lambda_0(\epsilon)
=\frac{2\sqrt{1+\epsilon^2}-2-\epsilon^2}4
=-\frac{\epsilon^4}{16}+O(\epsilon^6)
\quad\text{as}\quad\varepsilon\to0.
\end{equation}

Let us now test Theorem~\ref{main theorem}
by comparing the asymptotic formula from this theorem
with formula (\ref{eigenvalue quartic asymptotic formula}).
Substituting (\ref{metric quartic}) into (\ref{definition of tensor h}) we get
\[
h_{\alpha\beta}(x^1)=
\begin{pmatrix}
0&\cos x^1&\sin x^1\\
\cos x^1&0&0\\
\sin x^1&0&0
\end{pmatrix}.
\]
Application of the Fourier transform (\ref{Fourier transform}) gives us
\begin{equation}
\label{Example of quartic dependence h hat}
\hat h_{\alpha\beta}(m_1)=
\begin{cases}
\begin{pmatrix}
0&\frac12&-\frac i2\\
\frac12&0&0\\
-\frac i2&0&0
\end{pmatrix}
\quad&\text{for}\quad m_1=1,
\\
{}&{}
\\
\qquad\quad
0\quad&\text{for}\quad m_1=2,3,\ldots.
\end{cases}
\end{equation}
Substituting
(\ref{Example of quartic dependence h hat})
into
(\ref{formula for c sum axisymmetric})
we get
$c=0$, in agreement with  (\ref{eigenvalue quartic asymptotic formula}).

\section{The eta invariant}
\label{The eta invariant}

Let $H$ be a first order self-adjoint elliptic $m\times m$ matrix classical pseudodifferential
operator acting on $m$-columns of complex-valued half-densities over
a compact $n$-dimensional manifold $M$ without boundary.
Here ellipticity is understood as the nonvanishing of the
determinant of the principal symbol of $H$, see \cite{jst_part_a}.
The \emph{eta function} of
$H$ is defined as
\begin{equation}
\label{definition of eta function}
\eta_H(s):=\sum\frac{\operatorname{sign}\lambda}{|\lambda|^s}\,,
\end{equation}
where summation is carried out over all nonzero eigenvalues
$\lambda$ of $H$, and $s\in\mathbb{C}$ is the independent variable.
Asymptotic formulae for the counting function imply that
the series (\ref{definition of eta function}) converges absolutely for
$\operatorname{Re}s>n$ and defines a holomorphic function in this half-plane.
It is known
\cite{atiyah_part_3}
that the eta function extends meromorphically to the whole $s$-plane.
Moreover, it is known, see Theorem 4.5 in \cite{atiyah_part_3},
that if the dimension $n$ is odd, then the eta function
is holomorphic at $s=0$.
This justifies, for odd $n$,
the definition of the \emph{eta invariant} as the real number $\eta_H(0)$.
The eta invariant $\eta_H(0)$ is the traditional measure of spectral
asymmetry of the operator $H$.

If we have only a finite number of eigenvalues
(i.e.~if we are looking at an Hermitian matrix
rather than a differential operator) then
the eta invariant is an integer number:
it is the number of positive eigenvalues minus
the number of negative eigenvalues.
However, in the case of a differential operator there is no reason
for the eta invariant to be integer.
The basic example \cite{atiyah_short_paper} is that of the
scalar ordinary differential operator
$\,H(\epsilon):=-i\frac d{dx^1}+\epsilon\,$
acting on the unit circle parameterized
by the  cyclic coordinate $x^1$ of period $2\pi$,
with $\epsilon$ being a real parameter.
It is known \cite{atiyah_short_paper}
that the eta invariant $\eta_{H(\epsilon)}(0)$
of this ordinary differential operator
is the odd 1-periodic function defined by the
formula $\eta_{H(\epsilon)}(0)=1-2\epsilon$ for $\epsilon\in(0,1)$.
In particular, we have $\eta_{H(0)}(0)=0$ and
$\lim\limits_{\epsilon\to0^\pm}\eta_{H(\epsilon)}(0)=\pm1$.

The current state of affairs (from an analyst's perspective)
in the subject area of zeta/eta functions of elliptic operators is
described in detail in the two papers
\cite{grubb_and_seeley_1995,grubb_and_seeley_1996}.
Let us highlight a few facts.
\begin{itemize}
\item
The key results are Theorem~2.7 from \cite{grubb_and_seeley_1995}
and Proposition~2.9 from \cite{grubb_and_seeley_1996}. Arguing
along the lines of \cite{atiyah_part_3} one can recover from these
results, in a rigorous analytic fashion,
properties of the eta function.
\item
The eta function is holomorphic at $s=0$ in any dimension $n\in\mathbb{N}$
(i.e. without the assumption of $n$ being odd).
This fact was proved by P.~B.~Gilkey \cite{gilkey_1981}.
\item
The seminal paper of R.~T.~Seeley \cite{seeley}
contained a small mistake:
see page 482 in \cite{grubb_and_seeley_1995}
or Remark 2.6 on page 39 in \cite{grubb_and_seeley_1996}
for details.
\end{itemize}

The more recent survey papers
\cite{grubb_2005,grubb_2006}
provide an overview of the subject.

Let us denote our massless Dirac operator on half-densities by
$A(\epsilon)$, where $\epsilon\in\mathbb{R}$
is the small parameter appearing in our metric $g_{\alpha\beta}(x;\epsilon)$.
Theorem \ref{main theorem} implies the following corollary.

\begin{corollary}
\label{corollary about eta function}
Suppose that the coefficient $c$ defined by formula
(\ref{formula for c sum}) is nonzero. Then
\begin{equation}
\label{limit of eta function}
\lim_{\epsilon\to0}\eta_{A(\epsilon)}(0)=2\operatorname{sign}c\,.
\end{equation}
\end{corollary}

Note that we have $\eta_{A(0)}(0)=0$,
so formula (\ref{limit of eta function})
implies that the function $\eta_{A(\epsilon)}(0)$
is discontinuous at $\epsilon=0$.

\

\emph{Proof of Corollary \ref{corollary about eta function}\ }
Put
$
f(\epsilon,t):=
\operatorname{Tr}
\left[
A(\epsilon)\,
e^
{
-t(A(\epsilon))^2
}
\right]$,
where $t>0$ and $\,\operatorname{Tr}\,$ is the operator
(as opposed to pointwise) trace\footnote{
The paper \cite{bismut_and_freed} to which we are about
to refer to actually deals with pointwise estimates,
i.e.~the trace in \cite{bismut_and_freed} is understood as
the matrix trace of the integral kernel on the diagonal
at a given point $x$ of the manifold.
We do not need pointwise estimates for the proof of
Corollary~\ref{corollary about eta function}.
}.
Having fixed~$\epsilon$, let us examine the behaviour of
$f(\epsilon,t)$ as $t\to0^+$.
For a generic first order pseudodifferential operator $A(\epsilon)$ we would
have $f(\epsilon,t)=O(t^{-2})$.
However, as explained in Chapter II of \cite{bismut_and_freed},
the Dirac operator in odd dimensions is very special and there are
a lot of cancellations when one computes the asymptotic expansion for
$f(\epsilon,t)$ as $t\to0^+$.
Namely, it was shown in \cite{bismut_and_freed} that
\begin{itemize}
\item
$f(\epsilon,t)=O(\sqrt t\,)$ as $t\to0^+$,
\item
$\eta_{A(\epsilon)}(s)$ is holomorphic in the half-plane $\operatorname{Re}s>-2\,$, and
\item
\begin{equation}
\label{eta function via trace}
\eta_{A(\epsilon)}(s)=\frac1
{\Gamma\left(\frac{s+1}2\right)}
\int_0^{+\infty}t^{(s-1)/2}\,f(\epsilon,t)\,dt
\quad\text{for}\quad
\operatorname{Re}s>-2\,.
\end{equation}
\end{itemize}
See also Section 1 in \cite{woj1999}.

Formula (\ref{eta function via trace}) implies
\begin{equation}
\label{eta invariant via trace}
\eta_{A(\epsilon)}(0)=\frac1
{\sqrt\pi}
\int_0^{+\infty}\,\frac{f(\epsilon,t)}{\sqrt t}\ dt\,,
\end{equation}
so in order to prove Corollary~\ref{corollary about eta function}
we need to examine the behaviour of the integral
(\ref{eta invariant via trace}) as $\epsilon\to0$.

Let us denote by $\lambda_0(\epsilon)$
the eigenvalue of the operator $A(\epsilon)$ with smallest modulus
and by $E_0(\epsilon)$ the orthogonal projection onto the
corresponding 2-dimensional eigenspace.
Put
\[
A_0(\epsilon):=\lambda_0(\epsilon)\,E_0(\epsilon)\,,
\qquad
\tilde A(\epsilon):=A(\epsilon)-A_0(\epsilon),
\]
\[
f_0(\epsilon,t):=
\operatorname{Tr}
\left[
A_0(\epsilon)\,
e^
{
-t(A_0(\epsilon))^2
}
\right]
=
2\,\lambda_0(\epsilon)\,
e^
{
-t(\lambda_0(\epsilon))^2
},
\]
\[
\tilde f(\epsilon,t):=
\operatorname{Tr}
\left[
\tilde A(\epsilon)\,
e^
{
-t(\tilde A(\epsilon))^2
}
\right]=f(\epsilon,t)-f_0(\epsilon,t).
\]
Then formula (\ref{eta invariant via trace}) can be rewritten as
\begin{multline}
\label{3 terms}
\eta_{A(\epsilon)}(0)=
\frac1{\sqrt\pi}
\int_0^1\,\frac{f(\epsilon,t)}{\sqrt t}\ dt\,
+\,
\frac1{\sqrt\pi}
\int_1^{+\infty}\,\frac{\tilde f(\epsilon,t)}{\sqrt t}\ dt\,
\\
+\,
\frac2{\sqrt\pi}
\int_1^{+\infty}\,\frac{
\lambda_0(\epsilon)\,
e^
{-t(\lambda_0(\epsilon))^2}
}
{\sqrt t}\ dt\,.
\end{multline}
The three terms in the RHS of (\ref{3 terms}) are functions of the
parameter $\epsilon$ and we shall now examine how they depend on
$\epsilon$.

The first term in the RHS of (\ref{3 terms}) is continuous at
$\epsilon=0$ because asymptotic formulae for
$f(\epsilon,t)$ as $t\to0^+$ are uniform in $\epsilon$.
This follows from the construction of heat kernel type
asymptotics for $t\to0^+$: the algorithm
is straightforward and examination of this algorithm shows that the
asymptotic coefficients and remainder term depend on additional
parameters in a continuous fashion.

The second term in the RHS of (\ref{3 terms}) is continuous at
$\epsilon=0$ because the eigenvalues of the operator
$\tilde A(\epsilon)$ depend on $\epsilon$ continuously
and because all these eigenvalues,
bar one double eigenvalue, are uniformly separated from zero.
The double eigenvalue in question is identically zero as a function
of $\epsilon$ and does not contribute to the second term in the RHS of
(\ref{3 terms}).

Thus, the proof of formula (\ref{limit of eta function}) reduces to
the proof of the statement
\begin{equation}
\label{limit of eta function reduced}
\frac2{\sqrt\pi}
\,\lim_{\epsilon\to0}\,
\int_1^{+\infty}\,\frac{
\lambda_0(\epsilon)\,
e^
{-t(\lambda_0(\epsilon))^2}
}
{\sqrt t}\ dt\,
=2\operatorname{sign}c\,.
\end{equation}
But formula (\ref{limit of eta function reduced})
is an immediate consequence of formula (\ref{main theorem formula}).~$\square$

\

The geometric meaning of the eta invariant of the Dirac operator
acting over a compact oriented
Riemannian manifold of dimension $4k-1$, $k\in\mathbb{N}$,
has been extensively studied in
\cite{atiyah_short_paper,atiyah_part_1,atiyah_part_2,atiyah_part_3}.
We are, however, unaware of publications dealing specifically
with the Dirac operator on a 3-torus,
though certain 2-torus bundles over a circle were examined in
\cite{atiyah_annals}.

\section*{Acknowledgments}

The authors are grateful to M.~F.~Atiyah, P.~B.~Gilkey, G.~Grubb and J.~D.~Lotay
for advice on the eta invariant,
and to B.~Ammann and M.~Hortacsu
for drawing our attention to the papers
\cite{ammann2009,ammann2011, hortacsu}.

\end{document}